\documentclass[preprint,12pt]{elsarticle}

\newcommand{\beq}{\begin{equation}}
\newcommand{\eeq}{\end{equation}}
\newcommand{\bqy}{\setlength{\arraycolsep}{0.0em}\begin{eqnarray}}
\newcommand{\eqy}{\end{eqnarray}\setlength{\arraycolsep}{5pt} }
\newcommand{\pf}[2]{\frac{\partial #1}{\partial #2}}
\newcommand{\pfs}[2]{\frac{\partial^2 #1}{\partial #2^2}}

\newtheorem{thm}{Theorem}[section]
\newtheorem{defn}{Definition}[section]
\newtheorem{rem}{Remark}[section]

\newtheorem{cor}{Corollary}[section]
\newtheorem{lem}{Lemma}[section]
\newproof{prf}{Proof}

\usepackage{graphicx}
\usepackage{appendix}
\usepackage{arydshln}

\usepackage{amssymb}
\journal{System and Control Letters}

\begin{document}

\begin{frontmatter}


%\title{Global Adaptive Extremum Seeking Control for General Multidimensional Systems}
\title{ES-MRAC: A New Paradigm For Adaptive Control}
%\author[PU]{Poorya Haghi\corref{cor1}\fnref{fn1}}
\author[PU]{Poorya Haghi\corref{cor1}}
\ead{phaghi@purdue.edu}
%\fntext[fn1]{Poorya Haghi is with the School of Mechanical Engineering, Purdue University.}

%\author[PU]{Kartik B. Ariyur\fnref{fn2}}
\author[PU]{Kartik B. Ariyur}
\ead{kariyur@purdue.edu}
%\fntext[fn2]{Kartik B. Ariyur is with the School of Mechanical Engineering, Purdue University.}

\cortext[cor1]{Corresponding Author}
\address[PU]{School of Mechanical Engineering, 585 Purdue Mall, West Lafayette, IN, 47907}


\begin{abstract}

We develop a method for the model reference adaptive control (MRAC) of LTI systems via Extremum Seeking (ES). Our proof of global asymptotic tracking enables design of the adaptive controller to satisfy averaging requirements, and convergence of tracking error. Our method is novel, with the additional advantage that no perturbations need be added to the reference trajectory and no {\it a priori} knowledge of parameter signs is needed. We illustrate our results for a simulated second order system.

\end{abstract}


\end{frontmatter}

% \linenumbers

\section{Introduction}
\label{Intro}

Many engineered systems such as robots carrying objects with unknown inertias, or unmanned vehicles subject to uncertain forces, are yet expected to guarantee performance. Dealing with such systems has motivated the problems of  adaptive control. Adaptive control has a long history, dating back to the 1922 paper of Leblanc \cite{Leblanc}, whose scheme may have been the first ``adaptive'' controller. Designing autopoilots for aircraft motivated adaptive control in the 1950s \cite{Greg}, followed by developments of self-adjusting schemes such as M.I.T. rule \cite{MIT} and gradient estimation \cite{Kokotovic1964, KRP1985}. Over the years, several solutions have been provided to this fundamental problem \cite{ Astrom, Narendra, Egardt, Goodwin, Slotine, kkk, IK}.

Many methods of adaptive control have been proposed, e.g. extremum seeking \cite{Morosanov, Ostrovskii, kartik1}, self tuning regulators  \cite{Astrom, AW1973}, direct and indirect adaptations \cite{Ioannou}, and adaptive back-stepping \cite{kkk, KKM1991, KKM1991II}. Our motivation is control problems that require both transient and steady state performance. Therefore, we chose to adapt model reference control (MRC) since it can regulate the transients. MRC uses a reference model to specify the ideal response of the plant.
In this paper, we use extremum seeking (ES) loops, based on sinusoidal perturbations to optimize cost functions based on the tracking error of MRC. We prove global asymptotic tracking of the adaptive system, and develop systematic design guidelines to satisfy the conditions of the stability proof. This work complements our initial attempt at ES-MRAC \cite{ACCpaper}.
%The idea is to combine the methods of  ES and MRAC to create the so-called extremum seeking model reference adaptive control (ES-MRAC). We employ sinusoidal perturbations to perturb the output of a cost function. Extremum seeking is then used to optimize this cost function. This would lead to tracking of the reference or ideal model. 

%We first proposed this method in \cite{ACCpaper}, where special cases of first order systems. However, the proof was limited to averaged equations only. 
%Although the analysis is rather difficult, 
ES-MRAC confers advantages over classic MRAC: avoiding perturbation of the reference signal, and imposing no requirements on the signs of parameters. 
ES-MRAC is a novel paradigm for adaptive control and opens up many interesting theoretical and practical problems, which we list in our conclusions.

Our paper is organized as follows. We supply essential background from prior art in section \ref{basics}, our main result in section \ref{main}, and an illustrative example in section \ref{secondorder}. In section \ref{conclusions}, we list some of the questions opened up by our results.

\section{Essential Prerequisites}
\label{basics}

This section introduces some of the preliminary notions, definitions, and prior results required for thorough understanding of the paper. Our goal is to acquaint the reader with the tools that we use in section \ref{main}.

\subsection{Extremum Seeking}

Extremum seeking is a powerful tool for obtaining the extremum (minimum or maximum) value of a map. Hence, it is used in many control applications, where the reference-to-output map has an extremum value. 

\begin{figure}[btp]
\centering
\includegraphics[height = 4.0cm]{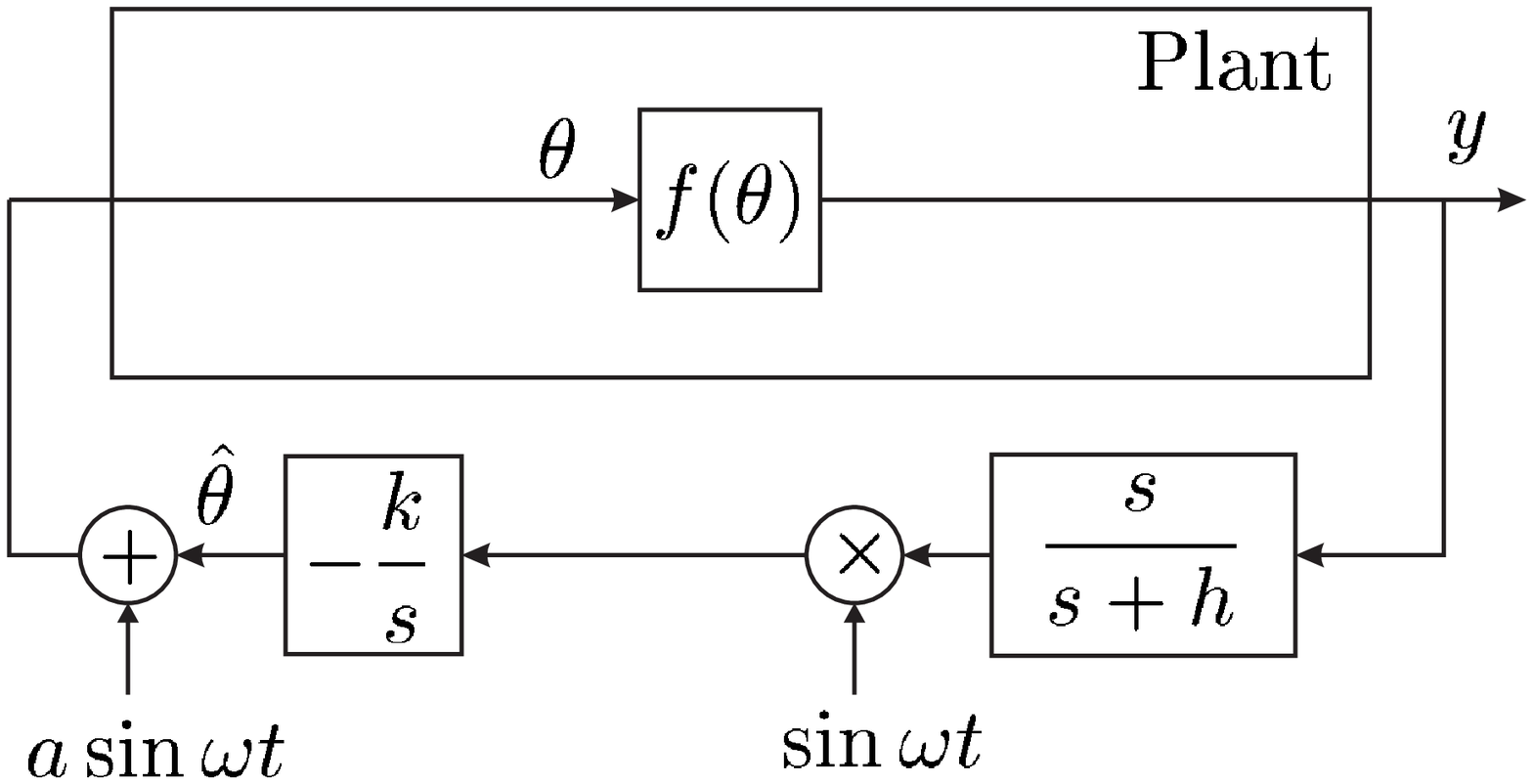}
\caption{The basic extremum seeking scheme with sinusoidal perturbations.}
\label{fig:BasicExtremumSeeking}
\end{figure}

Suppose we have a system with an input $\theta$, and an output $y$, which is an unknown function of the input, $y=f(\theta)$. Without loss of generality, we assume that the mapping, $f(\theta)$, has a minimum value. Extremum seeking deals with the problem of finding the optimal input, $\theta^*$, that would drive the output, $f(\theta)$, to its minimum value, $f^*$. The basic extremum seeking scheme is shown in Figure \ref{fig:BasicExtremumSeeking}. The perturbation signal $\sin\omega t$, provides a measure of gradient information of the map $f(\theta)$. The result is summarized as follows.

\begin{thm}[see \cite{kartik1}]
For the system in Figure \ref{fig:BasicExtremumSeeking}, the output error $y-f^*$ achieves local exponential convergence to an $O(a^2+1/\omega^2)$ neighborhood of the origin, provided the perturbation frequency $\omega$ is sufficiently large, and $\frac{1}{1+L(s)}$ is asymptotically stable, where $L(s)=\frac{kaf''}{2s}$.
\end{thm}

Our objective in this paper is to develop an adaptive controller that minimizes both the tracking errors and the errors in estimating the unknown parameters. Thus, we use extremum seeking as a means of optimizing a nonlinear cost function of the errors. In simple words, this cost function plays the role of the mapping $f(\theta)$ in Figure \ref{fig:BasicExtremumSeeking}. This process yields an adaptation law which  updates the gains of the controller. The use of a nonlinear cost function, along with the time dependent perturbation signals, lead to a non-autonomous, nonlinear set of equations. We use averaging and Lyapunov theory for stability analysis. Averaging removes time dependence from the equations and makes standard Lyapunov analysis possible. We reproduce essential definitions and theorems here for completeness and convenience. 

\begin{defn}[Averaging, \cite{Khalil}]
\label{def-averaging}
Consider the non-autonomous system $\dot x=\varepsilon f(t,x,\varepsilon)$, where $\varepsilon$ is a small positive constant, and $x\in D\subset R^n$. Suppose that $f$ is $T$-periodic in $t$, i.e. for all $t\geq 0$, $f(t+T,x,\varepsilon)=f(t,x,\varepsilon)$. We obtain the ``average system'' by \beq
\label{eq-averaging}
\dot{x}_{\mbox{\scriptsize{av}}}=\varepsilon f_{\mbox{\scriptsize{av}}}(x_{\mbox{\scriptsize{av}}}),
\eeq
where $f_{\mbox{\scriptsize{av}}}(x)=\frac{1}{T}\int_0^Tf(\tau,x, 0) d\tau$.
\end{defn}

%The average system defined above is autonomous and less complex, but obviously different from the original system. However, under certain conditions, the solution of the original system, $\dot x=\varepsilon f(x,t,\varepsilon)$, and the solution of the average system, $\dot x=\varepsilon f_{\mbox{\scriptsize{av}}}(x)$, differ no more than $O(\varepsilon)$. One theorem that provides such conditions, will be used later in the paper. 

%Once the autonomous equations  are derived, one can use various tools to study stability of the system. We use the direct method of Lyapunov and Barbalat's lemma. 

%\subsection{Stability}

\begin{thm}[see \cite{Khalil}]
\label{thm-averaging}
Consider the system 
\beq
\label{eq-app-1}
\dot x = \varepsilon f(t,x,\varepsilon),\qquad \varepsilon>0,
\eeq 
where $f$ and its partial derivatives with respect to $(x,\varepsilon)$ up to the second order are continuous and bounded for $(t,x,\varepsilon)\in [0,\infty)\times D_0\times [0,\varepsilon_0]$, for every compact set $D_0\subset D$, where $D\subset R^n$ is a domain. Suppose $f$ is $T$-periodic in $t$ for some $T>0$. Let $x_{\mbox{\scriptsize{av}}}(\varepsilon t)$ and $x(t,\varepsilon)$ denote the solutions of (\ref{eq-averaging}) and (\ref{eq-app-1}), respectively. If $x_{\mbox{\scriptsize{av}}}(\varepsilon t)\in D\ \forall \ t\in[0,\infty)$ and $x(0,\varepsilon)-x_{\mbox{\scriptsize{av}}}(0) = O(\varepsilon)$, then there exists $\varepsilon^*>0$  such that for all $0<\varepsilon<\varepsilon^*$, $x(t,\varepsilon)$ is defined and 
\beq
x(t,\varepsilon)-x_{\mbox{\scriptsize{av}}}(\varepsilon t)=O(\varepsilon)\quad\mbox{on}\quad [0,\infty)
\eeq 
\end{thm}

We assume that the reader is familiar with the standard theorems regarding global asymptotic stability. Hence, we shall not provide such theorems here (see section 4.1 of \cite{Khalil} for example). 
%To this end, we use the following theorem from \cite{Khalil}.
%\begin{thm}[see \cite{Khalil}]
%\label{thm-global-asymptotic-stability}
%Let $x=0$ be an equilibrium point for $\dot x = f(x)$. Let $V:{R}^n\to {R}$ be a continuously differentiable function such that 
%\bqy
%V(0) = 0\quad\mbox{and}\quad V(x)>0,\ \forall x \not = 0\\
%||x||\to\infty \quad\Rightarrow\quad V(x)\to\infty\\
%\dot V(x)<0,\quad \forall x \not = 0
%\eqy 
%then $x=0$ is globally asymptotically stable. 
%\end{thm}

We are now equipped with the tools required to carry out a rigorous analysis on the use of the method of extremum seeking in adaptive control. The problem formulation and the main results are described in the following section.

\section{Main Results}
\label{main}
%Observing the results of the first and second order systems, poses the question whether there exist a similar trend in higher order dimensions, and whether the proposed control method will yield stable equilibria. 

In this section, we generalize our results from \cite{KartikCurrentPaper} to a single input multi output (SIMO) LTI system of order $n$. We show that with the proposed adaptation procedure, an LTI system of an arbitrary order $n$ with full-state measurement, can achieve global asymptotic tracking for all states. 

Suppose that we have a linear system of order $n$, with control input $u\in {R}$. Without loss of generality, we assume that the governing equations are written in the form, 
\bqy
\label{eq-system}
a_ny^{(n)}+a_{n-1}y^{(n-1)}+\ldots+a_0y = u,
\eqy 
where $y^{(i)}\in {R},\ (i=0,\ldots,n)$ are measurable states, and $a_0,a_1,\ldots,a_n$ are system parameters. Furthermore, assume that all the parameters are unknown to the designer, i.e. $a_i,i=0,\ldots,n$ are unknown. Note that the usual MRAC requires the sign of $a_n$ to be known, whereas this methodology puts no restriction on the sign of $a_n$. Hence we can deal with situations in which no {\it a priori} knowledge on the sign of parameters is given. The objective is to design a control law  to track the reference model
\beq
\label{eq-reference-model}
a_{mn}y_m^{(n)}+a_{m(n-1)}y_m^{(n-1)}+\ldots+a_{m0}y_m=r(t),
\eeq 
where $a_{mn}p^{n}+a_{m(n-1)}p^{n-1}+\ldots+a_{m0}$ is Hurwitz. In the above equation, $r(t)$ denotes the reference signal, $y_m^{(i)}\in {R},\ (i=0,\ldots,n)$ are the states of the reference model, and $a_{m0},a_{m1},\ldots,a_{mn}$ are known constants.

The method that we propose, provides a control law, $u$, such that the system dynamics, (\ref{eq-system}), will follow the reference model ,(\ref{eq-reference-model}), and uses the same tracking error as MRAC. Uncertainties are then, taken into account by adaptation of control parameters via Extremum Seeking (ES) loops for each parameter. Hence, we call it ES-MRAC. %We propose the control law and the adaptation law as follows.
Prior to introducing the stability theorem, we provide the following definitions.

\begin{defn}
The `auxiliary signal', $z(t)$, is defined as follows
\beq
\label{eq-z-definition}
z(t) \triangleq y_m^{(n)}-\beta_{n-1}e^{(n-1)}-\ldots-\beta_0e
\eeq 
with the design paramters, $\beta_i$, chosen such that the polynomial $p^n+\beta_{n-1}p^{n-1}+\ldots+\beta_0$ is Hurwitz, and $y_m$ being the reference model as governed by (\ref{eq-reference-model}).
\end{defn}

\begin{figure}[btp]
\centering
\includegraphics[height = 3.5cm]{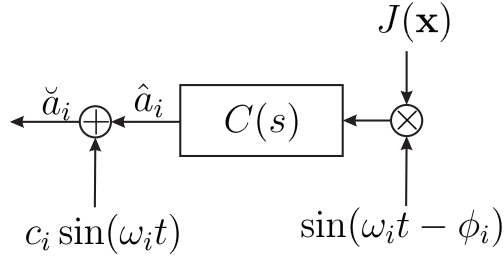}
\caption{ES block for updating parameters.}
\label{Fig-ESMRAC-General}
\end{figure}

The system parameter estimate is denoted by $\hat{a}_i$, ($i=0,\ldots,n$), and is adapted via an extremum seeking block as shown in Figure \ref{Fig-ESMRAC-General}, with a suitable choice of a cost function $J$, a compensator $C(s)$, and a probing frequency $\omega_i$.

\begin{defn}
The `perturbed estimate' of a parameter is denoted by $\breve{a}_i$, ($i=0,\ldots,n$), and defined  as the estimate of the parameter $a_i$ after being perturbed by the sinusoidal signal. In other words,
\beq
\label{eq5}
\breve{a}_i\triangleq\hat{a}_i+c_i\sin{\omega_i t}
\eeq
\end{defn}
\begin{defn}
The `control input' is defined as
\beq
\label{eq-ControlInput}
u \triangleq \breve{a}_nz(t)+\breve{a}_{n-1}y^{(n-1)}+\ldots+\breve{a}_0y
\eeq 
where $y$ is the system state as governed by (\ref{eq-system}), and $z(t)$ is the auxiliary signal.
\end{defn}
 
\begin{defn}
\label{def-trackingerror}
The `tracking error' is defined as the difference between the system state and the reference model state, i.e. 
\beq
\label{eq-tracking-error}
e\triangleq y-y_m
\eeq 
The `tracking error vector' refers to ${\bf x}=[x_1,\ldots,x_n]^T\triangleq\left[e,\dot e,\ldots,e^{(n-1)}\right]^T$.
\end{defn}

\begin{thm}
\label{thm-nth-order}
For an LTI system of order $n$ given by (\ref{eq-system}), with the control input $u$, (\ref{eq-ControlInput}), and the adaptation law for parameter estimate $\hat{a}_i$, $i=0,\ldots,n$, as given by the ES block shown in Figure \ref{Fig-ESMRAC-General},
let the cost function in Figure~\ref{Fig-ESMRAC-General} be 
\beq
\label{eq-J-definition}
J = \frac{1}{2}\left[{\bf q}^T{\bf x}\right]^2 = \frac{1}{2}\left[\sum_{k=1}^nq_ix_i\right]^2,
\eeq 
where ${\bf q}=[q_1,\ldots,q_n]^T$ is the vector of weighting factors, and ${\bf x}$ is the tracking error vector, and let the compensator in Figure~\ref{Fig-ESMRAC-General} be 
\beq
\label{eq-Compensator}
C(s) = -g_i\left(\frac{1+d_is}{s}\right)
\eeq 
Furthermore, let the probing frequency for each ES loop be given by $\omega_i=n_i\omega$, $n_i\in N$, $i=0,\ldots,n$ with $n_i\not=n_j$ for $i\not=j$. Then, there exist design parameters $c_i,d_i,g_i,\omega_i$, $\beta_j$, and $q_k$, ($i=0,\ldots,n$, $j=0,\ldots,n-1$, $k=1,\ldots,n$), such that this setup will guarantee global asymptotic convergence of the tracking error vector $\bf x$, to an $O(1/\omega)$ neighborhood of the origin, provided that the probing  frequency, $\omega$, is sufficiently large. 
\end{thm}

\begin{rem}
In the above theorem, design parameters refer to the following.
\begin{itemize}
\item weighting factors of the cost function ${\bf q}=[q_1,\ldots,q_n]^T$, 
\item compensator gains, $g_i$ and $d_i$,
\item $\beta_0,\ldots, \beta_{n-1}$ in Eq. (\ref{eq-z-definition}),
\item $c_i$ and $\omega_i$ (probing amplitude and frequency), in Eq. (\ref{eq5}).
\end{itemize}
\end{rem}
The proof of this theorem is postponed to Appendix\ref{appnedixA}. However, we present the following lemma which helps prove the theorem.
The proof of this lemma provides insight into how we can choose the design parameters. 

\begin{lem}
\label{lem-GlobalAsymptoticAverage}
Let $\bf x$ be the tracking error vector, and let the average tracking error, $\bf x_{av}$, be calculated as in Definition \ref{def-averaging}. If the conditions of Theorem \ref{thm-nth-order} are satisfied, then ${\bf x_{av}}\to 0$ globally and asymptotically.
\end{lem}
%The proof of this theorem falls into two parts. In the first part, we use averaging and prove that the averaged system is globally asymptotically stable. The proof of this lemma provides insight into how we can choose the design parameters. In the second part of the proof, we use the results of Teel et. al. \cite{04,05}, to show that the global asymptotic stability for the averaged system, indeed implies, semi-global practical asymptotic stability for the actual system. This part of the proof will be postponed to \ref{appnedixA}.

\begin{prf}%[Global Asymptotic Stability of Averaged Equations]
We start by deriving the governing dynamics, which includes the state tracking error dynamics, and the parameter estimation error dynamics.

\subsection*{Governing Dynamics:}
To find the state tracking error dynamics, we substitute (\ref{eq5}) and (\ref{eq-z-definition}) into (\ref{eq-ControlInput}), and the resultant into (\ref{eq-system}). Defining the state tracking error as in (\ref{eq-tracking-error}), and the parameter convergence error as
\beq
\tilde{a}_i=a_i-\hat{a}_i\qquad i=0,\ldots,n\ ,
\eeq
we rearrange the resulting equation to show that the state tracking error dynamics is given by
\bqy
\nonumber
a_n[e^{(n)}&{}+{}&\beta_{n-1}e^{(n-1)}+\ldots+\beta_0e]\\
\nonumber
&{}={}&\sum_{k=0}^{n-1}(c_k\sin\omega_kt-\tilde{a}_k)y^{(k)}+(c_n\sin\omega_nt-\tilde{a}_n)z(t)
\eqy 

Using Definition \ref{def-trackingerror}, we rewrite the state tracking error dynamics in state space form as
\bqy
\label{eq-error-dynamics-time}
\dot {\bf x} = {\bf Ax}+\frac{1}{a_n}{\bf b}{\bf v}^T[\mathcal{S}-\tilde{{\bf a}}]
\eqy 
where
\begin{eqnarray}
\label{eq-A-Hurwitz}
{\bf A} = \left[\begin{array}{ccccc}
0&1&0&\ldots&0\\
0&0&1&\ldots&0\\
\vdots&\vdots&\vdots&\ddots&\vdots\\
0&0&0&\ldots&1\\
-\beta_0&-\beta_1&-\beta_2&\ldots&-\beta_{n-1}\end{array}\right],\quad {\bf b}=\left[\begin{array}{c}
0\\
0\\
\vdots\\
0\\
1\end{array}\right],
\end{eqnarray}
$\tilde{\bf a} = [\tilde{a}_0,\ldots,\tilde{a}_n]^T$, and $\mathcal{S}$ denotes the vector of perturbation signals given by $\mathcal{S}=[c_0\sin\omega_0t,\ldots,c_n\sin\omega_nt]^T$. In addition, the vector $\bf v$ is the regression vector given by ${\bf v} = [y,\dot y, \ldots,y^{(n-1)},z]^T$.

According to Figure \ref{Fig-ESMRAC-General}, the parameter estimation error is governed by
\bqy
\dot{\tilde{a}}_i = g_i(1+d_is)[\sin(\omega_it-\phi_i)J],\quad i=0,\ldots,n
\eqy 
Applying $s$ to the term in brackets,
\bqy
\label{eq-adaptation-time}
\dot{\tilde{a}}_i = g_i\Big(\sin(\omega_it-\phi_i)+d_i\omega_i\cos(\omega_it-\phi_i)\Big)J
{+}g_id_i\sin(\omega_it-\phi_i)\dot J
\eqy 
where $J$ is given by (\ref{eq-J-definition}), and $i=0,\ldots,n$. 
In order to write the last equation in vector form, we define a gain matrix ${\bf G}=\mbox{diag}(g_i)\in R^{(n+1)\times (n+1)}, i=0,\ldots,n+1$, and perturbation vectors, ${\bf d_1}(t)$ and ${\bf d_2}(t)$ 
\bqy
\nonumber	
{\bf d_1}(t)=\left[\begin{array}{c}
\sin(\omega_0t-\phi_0)+d_0\omega_0\cos(\omega_0t-\phi_0)\\
\sin(\omega_1t-\phi_1)+d_1\omega_1\cos(\omega_1t-\phi_1)\\
\vdots\\
\sin(\omega_nt-\phi_n)+d_n\omega_n\cos(\omega_nt-\phi_n)
\end{array}\right],\ 
{\bf d_2}(t)=\left[\begin{array}{c}
d_0\sin(\omega_0t-\phi_0)\\
d_1\sin(\omega_1t-\phi_1)\\
\vdots\\
d_n\sin(\omega_nt-\phi_n)\\
\end{array}\right].
\eqy 
Thus, we write (\ref{eq-adaptation-time}) as
\beq
\label{eq-adaptation-time2}
{\bf \dot{\tilde{a}}} = {\bf Gd_1}J+{\bf Gd_2}\dot J
\eeq 
Equations (\ref{eq-error-dynamics-time}) and (\ref{eq-adaptation-time2}), constitute the governing dynamics. These equations are non-autonomous since time appears explicitly in them. %Therefore, we use averaging to convert the governing dynamics into a set of autonomous equations.

\subsection*{Averaging:}

Let the greatest common factor of all the probing frequencies, $\omega_i$,$i=0,\ldots,n$, be denoted by $\omega$. In other words, $\omega_i=n_i\omega$, where $n_i\in {N}$, $i=0,\ldots,n$. Furthermore, assume that $n_i\not= n_j$ for $i\not =j$. This will guarantee orthogonality of the sinusoidal perturbations. Let the design parameters be chosen such that 
\bqy 
\label{eq-greatest-common-factor}
\omega&{}\gg{}& 1\\
O(d_i\omega_i)&{}={}&1,\quad i=0,\ldots,n\\
O(g_i)&{}={}&1,\quad i=0,\ldots,n
\eqy

In order to be able to perform averaging, we scale the time as follows. Suppose that $\tau = \omega t$, and define $\varepsilon=\left(1/\omega\right)\ll 1$. By substituting $\tau= \omega t$ into (\ref{eq-error-dynamics-time}), and using Definition \ref{def-averaging}, one can show that the averaged equation for tracking error dynamics is given by
\bqy
\label{eq-systemI}
\left(\frac{d{\bf x}}{d\tau}\right)_{\mbox{\scriptsize{av}}} = \varepsilon\left[{\bf Ax}_{\mbox{\scriptsize{av}}}-\frac{1}{a_n}{\bf b}{\bf v}^T_{\mbox{\scriptsize{av}}}\tilde{{\bf a}}_{\mbox{\scriptsize{av}}}\right]
\eqy 

Similarly, performing the same procedure on (\ref{eq-adaptation-time2}), one can show that the averaged equation for parameter estimation errors is given by
\beq
\label{eq-systemII}
\left(\frac{d\tilde{\bf a}}{d\tau}\right)_{\mbox{\scriptsize{av}}} = \varepsilon\frac{q_n}{2a_n}\mathcal{C}{\bf v}_{\mbox{\scriptsize{av}}}{\bf q}^T{\bf x}_{\mbox{\scriptsize{av}}},
\eeq 
where $\mathcal{C}=\mbox{diag}(g_id_ic_i\cos\phi_i)\in{R}^{(n+1)\times(n+1)}$, $i=0,\ldots,n$. Appendix\ref{appendixB} provides the details on how averaging can lead us from (\ref{eq-adaptation-time2}) to (\ref{eq-systemII}).

Note that the perturbation amplitude, $c_i$, is chosen so as to produce a measurable variation in the plant output at the corresponding frequency.

\subsection*{Lyapunov analysis:}

Now, we are ready to use Lyapunov stability analysis to prove convergence of averaged tracking error, ${\bf x}_{\mbox{\scriptsize{av}}}$, to zero.
Consider the following Lyapunov function
\beq
V = {\bf x}_{\mbox{\scriptsize{av}}}^T{\bf P}{\bf x}_{\mbox{\scriptsize{av}}}+2{\bf \tilde{a}}_{\mbox{\scriptsize{av}}}^T{\bf \Gamma}{\bf \tilde{a}}_{\mbox{\scriptsize{av}}},
\eeq 
where $\bf P$, and $\bf \Gamma$ are symmetric positive definite matrices.
Without loss of generality, we assume that $\varepsilon =1$, and conduct the stability analysis. We shall use the prime symbol to denote differentiation with respect to $\tau$. Taking the derivative of $V$ with respect to $\tau$ and noting that $\bf P$ and $\bf \Gamma$ are symmetric matrices, we can write
\beq
V'={\bf x}_{\mbox{\scriptsize{av}}}'^T{\bf P}{\bf x}_{\mbox{\scriptsize{av}}}+{\bf x}_{\mbox{\scriptsize{av}}}^T{\bf P}{\bf x}_{\mbox{\scriptsize{av}}}'+4{\bf \tilde{a}}_{\mbox{\scriptsize{av}}}^T{\bf \Gamma}{\bf \tilde{a}}_{\mbox{\scriptsize{av}}}'
\eeq
Substituting (\ref{eq-systemI}) and (\ref{eq-systemII}), and simplifying, we get
\bqy 
\nonumber
V' = {\bf x}_{\mbox{\scriptsize{av}}}^T({\bf A}^T{\bf P}+{\bf P}{\bf A}){\bf x}_{\mbox{\scriptsize{av}}}+\frac{2}{a_n}{\bf \tilde{a}}_{\mbox{\scriptsize{av}}}^T\left[q_n{\bf \Gamma}\mathcal{C}{\bf v}_{\mbox{\scriptsize{av}}}{\bf q}^T-{\bf v}_{\mbox{\scriptsize{av}}}{\bf b}^T{\bf P}\right]{\bf x}_{\mbox{\scriptsize{av}}}
\eqy 
which we write as
\bqy 
\label{eq17}
V' = -{\bf x}_{\mbox{\scriptsize{av}}}^TQ{\bf x}_{\mbox{\scriptsize{av}}}+\frac{2}{a_n}{\bf \tilde{a}}_{\mbox{\scriptsize{av}}}^T\left[q_n{\bf \Gamma}\mathcal{C}{\bf v}_{\mbox{\scriptsize{av}}}{\bf q}^T-{\bf v}_{\mbox{\scriptsize{av}}}{\bf b}^T{\bf P}\right]{\bf x}_{\mbox{\scriptsize{av}}}
\eqy 
with $Q = -({\bf A}^T{\bf P}+{\bf P}{\bf A})$  a positive definite matrix.
By choosing the parameters in $\bf P$, $\mathcal{C}$, $\bf q$, and $\bf \Gamma$ such that 
\beq
\label{eq-MatrixCondition}
q_n{\bf \Gamma}\mathcal{C}{\bf v}_{\mbox{\scriptsize{av}}}{\bf q}^T-{\bf v}_{\mbox{\scriptsize{av}}}{\bf b}^T{\bf P}=0,\quad \forall {\bf v}_{\mbox{\scriptsize{av}}}\in  R^{(n+1)}
\eeq 
%\beq 
%\label{eq18}
%q_n{\bf \Gamma}\mathcal{C}{\bf v}{\bf q}^T-{\bf v}{\bf b}^T{\bf P}=0
%\eeq
we get
\bqy
\label{eq-barbalat}
V'= -{\bf x}_{\mbox{\scriptsize{av}}}^TQ{\bf x}_{\mbox{\scriptsize{av}}}\leq 0,
\eqy 
that is, $V'$ is negative semi-definite. This implies that $V(\tau)\leq V(0)$, and therefore, $\bf x_{\mbox{\scriptsize{av}}}$, and $\bf \tilde{a}_{\mbox{\scriptsize{av}}}$ are bounded. Moreover, $\bf v_{\mbox{\scriptsize{av}}}$ is also bounded by definition, since all its components are linear combinations of elements of $\bf x_{\mbox{\scriptsize{av}}}$ and the reference model. Therefore, $V''=-2{\bf x}_{\mbox{\scriptsize{av}}}^T{\bf QA}{\bf x}_{\mbox{\scriptsize{av}}}+\frac{2}{a_n}{\bf \tilde{a}}_{\mbox{\scriptsize{av}}}^T{\bf v}_{\mbox{\scriptsize{av}}}{\bf b}^T{\bf Q}{\bf x}_{\mbox{\scriptsize{av}}}$ is also bounded. Hence $V'$ is uniformly continuous. Therefore, by Barbalat's lemma, $V'\to 0$ as $t\to\infty$. Hence, according to (\ref{eq-barbalat}), ${\bf x}_{\mbox{\scriptsize{av}}}\to 0$ as $t\to\infty$, i.e. tracking error and all its derivatives converge asymptotically to zero. Furthermore, since $V\to\infty$ as $||\bf x_{\mbox{\scriptsize{av}}}||\to\infty$, global asymptotic tracking is achieved. %according to theorem \ref{thm-global-asymptotic-stability}.  
\end{prf}

The following corollary, which follows from the above proof, provides some guidelines as how to choose the design parameters. 

\begin{cor}[Design]
\label{cor3-1}
 The design parameters $c_i,d_i,g_i,\omega_i$, $\beta_j$, and $q_k$, ($i=0,\ldots,n$, $j=0,\ldots,n-1$, $k=1,\ldots,n$) in Theorem \ref{thm-nth-order}, must be chosen such that the following holds
\begin{enumerate}
\item $\omega{}\gg{} 1$
\item \label{cond2} $O(d_i\omega_i){}={}1,\quad i=0,\ldots,n$
\item \label{cond3} $O(g_i){}={}1,\quad i=0,\ldots,n$
\item \label{cond4} The matrix ${\bf P}=[P_{ij}]$ is found by solving the identity ${\bf PA}+{\bf A}^T{\bf P}=-{\bf Q}$ for some positive definite matrix $\bf Q$, with $\bf A$ defined as in (\ref{eq-A-Hurwitz}).
\item Eigenvalues of the matrix $q_n\Gamma\mathcal{C}\in R^{(n+1)\times (n+1)}$ must satisfy
\beq
\label{eq-cor}
\lambda_1=\lambda_2=\ldots=\lambda_{(n+1)}=\frac{q_1P_{n1}+q_2P_{n2}+\ldots+q_nP_{nn}}{q_1^2+q_2^2+\ldots+q_n^2 }
\eeq 
For the special case where $\Gamma$ is diagonal, $\Gamma=\mbox{diag}(\gamma_0,\ldots,\gamma_n)$, this equation simplifies to
\beq
\label{eq:design}
q_ng_id_ic_i\cos\phi_i\gamma_i=\frac{q_1P_{n1}+q_2P_{n2}+\ldots+q_nP_{nn}}{q_1^2+q_2^2+\ldots+q_n^2 },\quad i=0,\ldots,n
\eeq 
\end{enumerate}
\end{cor}
For example, we start by choosing the weighting factors $q_i$, and a positive definite matrix $\bf Q$, and solve for $\bf P$ using condition \ref{cond4} above. Next, we choose large enough $c_i$'s to produce measurable variations in the plant output in the presence of noise. Then, for each parameter estimate, the gains are tuned by choosing $g_i$, and $d_i$ such that their orders of magnitude satisfy conditions \ref{cond2} and \ref{cond3}, and their values satisfy (\ref{eq:design}) for some positive definite diagonal matrix $\Gamma$.

%\begin{rem}
%The conditions in the above Corollary are based on the use of Averaging in the proof of Theorem \ref{thm-nth-order}. Therefore, one might be able to find a more relaxed  set of conditions, using a different method for proving global asymptotic stability.
%\end{rem}

\begin{rem}
The last condition in Corollary \ref{cor3-1}, follows from simplification and mathematical interpretation of (\ref{eq-MatrixCondition}) as follows. Right multiplying (\ref{eq-MatrixCondition}) by $\bf q$ and noting that ${\bf q}^T{\bf q}$ and ${\bf b}^T{\bf Pq}$ are scalars, we get
\beq
\label{eq-Amatrix}
\left(\Gamma\mathcal{C}-\frac{1}{q_n}\frac{{\bf b}^T{\bf Pq}}{{\bf q}^T{\bf q}}{\bf I}_{(n+1)}\right){\bf v_{av}}=0,\quad \forall {\bf v_{av}}\in  R^{(n+1)}
\eeq 
which is basically the characteristic equation of the matrix $\Gamma \mathcal{C}$, i.e. its eigenvalues are all identically equal to $\frac{1}{q_n}\frac{{\bf b}^T{\bf Pq}}{{\bf q}^T{\bf q}}$. %Since the last equation holds for all ${\bf v_{av}}\in R^{(n+1)}$, the null space of $\mathcal{A}$ must cover the space $R^{(n+1)}$. Hence, rank$(\mathcal{A})=0$, meaning that all eigenvalues of $\mathcal{A}$ must be zero. Equation (\ref{eq-cor}), is then, a direct result from the definition of $\mathcal{A}$ as given by (\ref{eq-Amatrix}).
\end{rem}

\begin{rem}[Geometrical Interpretation]
\label{rem-geometrical}
Equation (\ref{eq-cor}) can be interpreted as follows. Consider the vector of weighting factors $\bf q$, and the vector ${\bf P}^T{\bf b}$ in the hyper space $R^n$. Let the projection of vector ${\bf P}^T{\bf b}$ onto the vector $\bf q$ be denoted by $\bf P_1$. Since $\bf P_1$ is in the same direction as $\bf q$, we can write it as a coefficient of $\bf q$, as follows
\[{\bf P_1}=\kappa {\bf q}\quad \mbox{for some }\ \kappa\in R\]
This coefficient is given by
\[\kappa  = \frac{\langle{\bf P}^T{\bf b},{\bf q}\rangle}{\langle{\bf q},{\bf q}\rangle}=\frac{{\bf b}^T{\bf Pq}}{{\bf q}^T{\bf q}}\]
which, according to (\ref{eq-cor}), equals the eigenvalues of the matrix $q_n\Gamma\mathcal{C}$. Therefore, for a chosen diagonal matrix $\Gamma$, an increase in the projection length, $\kappa$, increases the eigenvalues of $\mathcal C$, hence an increase in the rate of change of parameter estimates (as seen by equation (\ref{eq-systemII})).
\end{rem}

\begin{cor}
If the reference signal $r(t)$ is a unit step, then ES-MRAC guarantees that 
\begin{itemize}
\item[(i)] $\forall j\in\{0,\ldots,n\}\quad\tilde{a}_{\mbox{\scriptsize{av}},j}\to k_j,\quad$ as $t\to\infty\quad$ for some $k_j\in R$, and
\item[(ii)] $k_0=0$.
\end{itemize} 
\end{cor}
\begin{prf}
In Lemma \ref{lem-GlobalAsymptoticAverage}, we showed that ${\bf x}_{\mbox{\scriptsize{av}}}\to 0 $ globally and asymptotically. Substituting ${\bf x}_{\mbox{\scriptsize{av}}}\to 0 $ into (\ref{eq-systemII}) proves the first part of the corollary. To prove the second part, we consider the steady state solution of (\ref{eq-systemI}). Substituting ${\bf x}_{\mbox{\scriptsize{av}}}\to 0 $ into (\ref{eq-systemI}) yields the following 
\beq
{\bf b}{\bf v}_{\mbox{\scriptsize{av}}}^T{\bf \tilde{a}}_{\mbox{\scriptsize{av}}}\to 0
\eeq
Substituting for ${\bf b}$, and ${\bf v}_{\mbox{\scriptsize{av}}}$, we get
\beq 
\label{eq-bva}
{\bf b}{\bf v}_{\mbox{\scriptsize{av}}}^T{\bf \tilde{a}}_{\mbox{\scriptsize{av}}} = \left[\begin{array}{cccc}
0&\ldots&0&0\\
\vdots&\ddots&\vdots&\vdots\\
0&\ldots&0&0\\
y_{\mbox{\scriptsize{av}}}&\ldots&y_{\mbox{\scriptsize{av}}}^{(n-1)}&z_{\mbox{\scriptsize{av}}} 
\end{array}\right]\left[\begin{array}{c}
\tilde{a}_{\mbox{\scriptsize{av}},0}\\
\vdots\\
\tilde{a}_{\mbox{\scriptsize{av}},n-1}\\
\tilde{a}_{\mbox{\scriptsize{av}},n}
\end{array}\right]\to 0
\eeq 
\end{prf}
Since ${\bf x}_{\mbox{\scriptsize{av}}}\to 0$, therefore, $y_{\mbox{\scriptsize{av}}}\to y_m$ as $t\to\infty$. Using final value theorem, one can easily show that for a unit step reference signal $r(t)$, $\lim_{t\to\infty}y_m(t)=1/a_{m0}$. Thus, (\ref{eq-bva}) will simplify to
\beq
\frac{1}{a_{m0}}\tilde{a}_{\mbox{\scriptsize{av}},0}\to 0
\eeq 
which is the same as $k_0=0$.

\section{Second Order Example}
\label{secondorder}

As a special case of higher order systems, consider the equation of motion for a linear second order system
\beq 
a_2\ddot y + a_1\dot y+a_0y\ = u
\eeq 
with unknown $a_0$, $a_1$ and $a_2$.
The control objective is to design a state feedback control law $u$, such that the system follows the reference dynamics given by
\beq 
a_{m2}\ddot{y}_m+a_{m1}\dot{y}_m+a_{m0}y_m = r(t)
\eeq 
We use theorem \ref{thm-nth-order} and Corollary \ref{cor3-1} to accomplish this. These results simplify the design process into a few easily carried out steps. In general, we divide the design process into two parts. That is, determination of control/adaptation schemes, and gain tuning.

\subsection{Determination of Control and Adaptation}
\subsubsection*{Step 1: Control}
Use $n=2$ in Theorem \ref{thm-nth-order} to provide the control law as follows
\beq
u = \breve{a}_2z+\breve{a}_1\dot y + \breve{a}_0 y
\eeq 
with the auxiliary signal
$
z(t) = \ddot{y}_m-\beta_1\dot e - \beta_0e
$
and $\breve{a}_n=\hat{a}_n+c_n\sin\omega_nt$, $n=0,1,2$.
\subsubsection*{Step 2: Adaptation}
For each of the unknown parameters, $a_0$, $a_1$, and $a_2$, set up the adaptation block as given by Figure \ref{Fig-ESMRAC-General}, using distinct frequencies. For each block, the cost function $J$ is given by (\ref{eq-J-definition}). 
\beq
J = \frac{1}{2}\left[{\bf q}^T{\bf x}\right]^2 = \frac{1}{2}\left[q_1e+q_2\dot e\right]^2,
\eeq 
and the compensator $C(s)$ is given by (\ref{eq-Compensator}).

Given numerical values, one can already implement the above laws, but it is very difficult to design proper gains such that the method performs well. Therefore, in what follows, we shall provide some insights into the gain tuning process.
\subsection{Gain Tuning}
In general, tuning the gains for optimal performance is very difficult in adaptive control, especially in MRAC. One advantage of ES-MRAC is that it provides some guidelines on gain tuning. In what follows, we shall use the results of Corollary \ref{cor3-1} in a step by step manner, to help us in the process of gain tuning. Note that these steps are only provided here as guidelines, and one can use a different procedure to tune the gains. Moreover, these steps help the designer with a quick estimate of the orders of magnitudes of different gains, but will not necessarily provide optimal gains. 

\subsubsection*{Step 1:}
Choose $\beta_0$ and $\beta_1$ such that $p^2+\beta_1p+\beta_0$ is Hurwitz. These parameters effect the rate of convergence of tracking error. Then, the $\bf A$ matrix is determined from (\ref{eq-A-Hurwitz})
\begin{eqnarray}
{\bf A} = \left[\begin{array}{cc}
0&1\\
-\beta_0&-\beta_1\end{array}\right],
\end{eqnarray}
\subsubsection*{Step 2:}
Choose the coefficients $q_1$ and $q_2$ of the cost function $J$, as defined in~(\ref{eq-J-definition}). 
\subsubsection*{Step3:}
Pick a positive definite matrix $\bf Q$ and solve the identity ${\bf PA}+{\bf A}^T{\bf P}=-{\bf Q}$  for the $2\times 2$ matrix $\bf P$.
\subsubsection*{Step 4:}
Choose the probing frequencies of the sinusoidal perturbations $\omega_0$, $\omega_1$, and $\omega_2$, such that the sinusoidal terms are orthogonal. Once these values are determined, pick the compensator gains $d_0$, $d_1$, and $d_2$ such that $O(d_i\omega_i)=1$, $i=0,1,2$.
\subsubsection*{Step 5:}
Finally, decide the values of $g_i$, $c_i$, $\phi_i$, and $\gamma_i$ such that (\ref{eq:design}) holds. Each of these design parameters contribute differently to the control problem. It is desirable to make $\gamma_i$'s small, since they have an inverse effect on the rate of convergence of parameters (See Remark \ref{rem-geometrical}). The values for $c_i$'s provide the amplitude of the perturbation signals. The $c_i$'s should be large enough to produce measureable variations in the plant output, but cannot be too large, since they can cause instability, and excitation of higher dynamics which is undesirable. Finally, the compensator gains, $g_i$'s, need to be large, since they can increase the rate of change of parameters, as seen by (\ref{eq-systemII}).

%%------------------------end of second order systems--------------------%%

\subsection{Numerical Example and Simulations}
\label{Simu}

We use the following numerical values for our simulations.
Let the true parameters of the system $a_2\ddot y+a_1\dot y + a_0y=u$, have the values $a_2=1$, $a_1=3$, and $a_0=6.25$. Let the model reference be given by $\ddot{y}_m+4.2\dot{y}_m + 9y_m=r$, which corresponds to a system with a natural frequency of $\omega_n = 3$ rad/sec and damping ratio of $\zeta = 0.7$. The initial conditions are assumed to be zero for the model reference and $[y(0), \dot y(0)]^T=[-0.1,0.2]^T$ for the plant. The reference input $r$ is assumed to be a unit step function.

\begin{figure}[btp]
\centering
\includegraphics[height = 6.8cm]{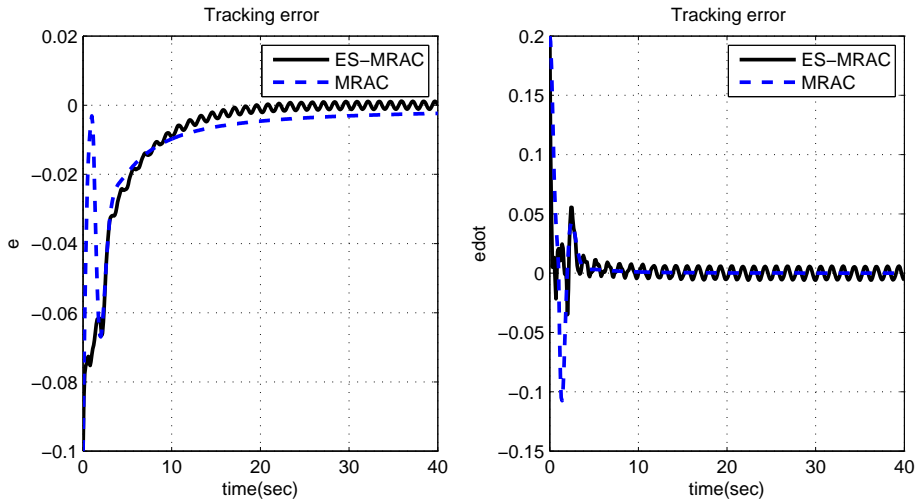}
\caption{Convergence of tracking error to zero for ES-MRAC in a 2nd order system.}
\label{fig:simulaiton1}
\end{figure}

\begin{figure}[btp]
\centering
\includegraphics[height = 6.8cm]{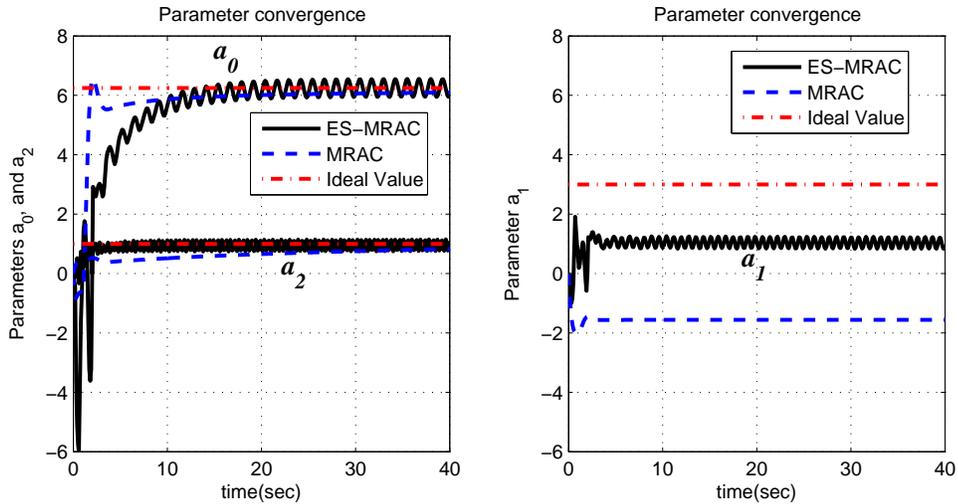}
\caption{Convergence of parameters. $a_{0}$, and $a_2$ reach their true values, while $a_1$ only reaches a constant.}
\label{fig:simulation2}
\end{figure}

It is assumed that there is no {\it a priori} knowledge of the ideal values of parameters. Figures \ref{fig:simulaiton1} and \ref{fig:simulation2} demonstrate  the performance of the system when the design parameters are chosen as follows: Let $\beta_0=9$, $\beta_1=3$, perturbation amplitudes $c_0=0.3$, $c_1=c_2=0.2$ perturbation frequencies $\omega_0=5$ rad/sec, $\omega_1=8$ rad/sec, $\omega_2=14$ rad/sec damping coefficients $d_1=d_2=d_3=0.1$, and gains $g_1=9000$, $g_2=3200$, $g_3=2000$. The matrix $\Gamma$ is chosen to be $\Gamma = 0.01{\bf I_{3}}$.

We compare our results with those of MRAC applied to the same system. Definition of $\Gamma$ in chapter 8 of \cite{Slotine} for MRAC, is the inverse of what we have defined. Therefore, we use $\Gamma = \mbox{diag}(90,90,50)$ for MRAC. Although this is not the exact inverse of the matrix $0.01{\bf I_{3}}$, it produces better results for MRAC.

\section{Conclusion}
\label{conclusions}
\noindent
We proposed
 a new approach to
 Model Reference Adaptive Control (MRAC), in which adaptation is carried through Extremum Seeking (ES), and we called it ES-MRAC. The proposed scheme was presented for a general class of LTI systems and proof of global asymptotic tracking was provided. The results of this paper open up many questions of practical relevance. %We list a few as follows. As an illustration, the adaptive control of a second order system was simulated using the proposed method. Results of the simulation, conformed to the theory we developed. 
Extensions to output feedback, direct adaptive control, and feedback linearizable systems may all be possible. Adaptive back-stepping for nonlinear systems and standard gradient and least square based adaptive controls may also have ES-MRAC counterparts. 
% Moreover, 

%% References


\bibliographystyle{elsarticle-num}
\bibliography{IEEEabrv,ExtremumBibFileAbb}






% \end{thebibliography}

\appendix
\renewcommand{\appendixname}{}

\section{Proof of Theorem \ref{thm-nth-order}} 
\label{appnedixA}

%Uniform boundedness is defined as follows as given by \cite{Khalil}.
%\begin{thm}
%Suppose that 
%\beq
%\label{eq-app-5}
%\dot x = f(t,x)
%\eeq
%where $f:[0,\infty)\times D\to \mathbb{R}^n$ is piecewise continuous in $t$ and locally Lipschitz in $x$ on $[0,\infty)\times D$, and $D\subset \mathbb{R}^n$ is a domain that contains the origin.
%
%The solutions of (\ref{eq-app-5}) are uniformly bounded if there exists a positive constant $c$, independent of $t_0\geq 0$, and for every $a\in (0,c)$, there is $\beta=\beta(a)>0$, independent of $t_0$, such that 
%\beq
%||x(t_0)||\leq a\quad\Rightarrow\quad ||x(t)||\leq \beta, \quad\forall t\geq t_0
%\eeq 
%\end{thm}

\begin{prf}%[Theorem \ref{thm-nth-order}]
%As shown in Lemma \ref{lem-GlobalAsymptoticAverage} through the Lyapunov analysis and Barbalat's lemma, the averaged system (\ref{eq-systemI}) and (\ref{eq-systemII}) is globally asymptotically stable in $\bf x_{\mbox{\scriptsize{av}}}$.

In order to study the properties of the actual system, we note that we can write the governing equations of the system, as $\dot x = f(x,t)$. That is, (\ref{eq-error-dynamics-time}) and (\ref{eq-adaptation-time2}) can be written as $\dot x = f(x,t)$, with $x = [{\bf x}^T,{\bf \tilde a}^T]^T$, and
\bqy
\label{eq-f-def}
f(t,x) = \left[\begin{array}{c}
\dot{\bf x}\\
\hdashline
\dot{\bf \tilde{a}}
\end{array}\right] = \left[\begin{array}{c}
{\bf Ax}+\frac{1}{a_n}{\bf b}{\bf v}^T[\mathcal{S}-\tilde{{\bf a}}]\\
\hdashline
{\bf Gd_1}J+{\bf Gd_2}\dot J
\end{array}\right]
\eqy 
As shown in section \ref{main}, we scale the time using $\tau = \omega t$, with $\omega \gg 1$, and $\varepsilon = 1/\omega \ll 1$. Thus we get
\beq
\frac{dx}{d\tau}=\varepsilon f(\tau,x)
\eeq 
As we can see from (\ref{eq-f-def}), $f$, $\pf{f}{x}$, and $\pfs{f}{x}$ are continuous. 
Furthermore, we note that on any compact set $D_0\subset D\subset R^{2n+1}$, $x\in D_0$ is bounded. That means that ${\bf x}$ and $\bf \tilde{a}$ are bounded. Therefore, according to (\ref{eq-f-def}), $f(t,x)$ (similarly $f(\tau,x)$) is bounded. Moreover, $f(\tau,x)$ as defined above is $2\pi$-periodic, since we are only using sinusoidal perturbations, and have defined $\omega$ to be the greatest common factor of all probing frequencies (see section \ref{main}, equation (\ref{eq-greatest-common-factor})).

According to Lemma \ref{lem-GlobalAsymptoticAverage}, for the averaged system (\ref{eq-systemI}) and (\ref{eq-systemII}), ${\bf x}_{\mbox{\scriptsize{av}}}\to 0$ globally and asymptotically, and ${\bf \tilde{a}}_{\mbox{\scriptsize{av}}}$ is bounded. Therefore, denoting the averaged solution for the overall system by $x_{\mbox{\scriptsize{av}}}(\varepsilon\tau) = [{\bf x}_{\mbox{\scriptsize{av}}}^T,{\bf \tilde a}_{\mbox{\scriptsize{av}}}^T]^T$, we see that the averaged system is globally bounded, i.e.  $x_{\mbox{\scriptsize{av}}}(\varepsilon\tau)\in D \ \forall \tau\in [0,\infty)$.

Hence according to Theorem \ref{thm-averaging}, solving the average system with the same initial state as the original system yields
\beq
\label{eq-app-2}
x(\tau,\varepsilon)-x_{\mbox{\scriptsize{av}}}(\varepsilon\tau)= O\left(\frac{1}{\omega}\right) \quad\mbox{on}\quad [0,\infty)
\eeq 
where we substituted $\varepsilon=1/\omega$. Equation (\ref{eq-app-2}) means that the solution to the nonautonomous system is only $O(1/\omega)$ different than the solution of the average system, i.e. ${\bf x}(\tau)-{\bf x}_{\mbox{\scriptsize{av}}}(\varepsilon\tau)=O(1/\omega)$, and ${\bf \tilde{a}}_{\mbox{\scriptsize{av}}}(\tau)-{\bf \tilde{a}}(\varepsilon\tau)=O(1/\omega)$. Since we showed in Lemma \ref{lem-GlobalAsymptoticAverage} that ${\bf x}_{\mbox{\scriptsize{av}}}\to 0$ globally and asymptotically, thus the vector ${\bf x}$ converges globally and asymptotically to on $O(1/\omega)$ neighborhood of the origin.	
\end{prf}

% ---------------- Details of Averaging the Adaptation Law------------

\section{Details of Averaging the Adaptation Law}
\label{appendixB}

As an example of how averaging is performed, we shall provide details of the averaging for the adaptation law. We will show how equation (\ref{eq-systemII}) is found by averaging (\ref{eq-adaptation-time2}).

We start by defining an averaging operator, using Definition \ref{def-averaging}.
\begin{defn}
The averaging operator, {\bf AVG}(.), for a $T$-periodic function $f(x,t,\epsilon)$ is defined as
\beq
\mbox{\bf AVG}\left(f(x,t,\epsilon)\right) = \frac{1}{T}\int_0^T f(x,\tau,0)d\tau
\eeq
\end{defn} 
\begin{rem}
The {\bf AVG} operator is linear
\[{\bf AVG}(c_1f_1+c_2f_2)=c_1{\bf AVG}(f_1)+c_2{\bf AVG}(f_2),\ \forall c_1,c_2\in R\]
\end{rem}

Next, we substite $J=\frac{1}{2}\left({\bf q}^T{\bf x}\right)^2$, and its derivative into (\ref{eq-adaptation-time2})
\bqy
{\bf \dot{\tilde{a}}} = \frac{1}{2}{\bf Gd_1}\left({\bf q}^T{\bf x}\right)^2+{\bf Gd_2}\left({\bf q}^T{\bf x}\right)\left({\bf q}^T{\bf \dot x}\right)
\eqy
Substituting for $\bf \dot x$ from (\ref{eq-error-dynamics-time}) would then give
\bqy
\label{eq-b2}
{\bf \dot{\tilde{a}}} = \frac{1}{2}{\bf Gd_1}\left({\bf q}^T{\bf x}\right)^2+{\bf Gd_2}{\bf q}^T{\bf x}{\bf q}^T\left({\bf Ax}+\frac{1}{a_n}{\bf b}{\bf v}^T[\mathcal{S}-\tilde{{\bf a}}]\right)
\eqy

Recall that we defined the perturbation frequencies as $\omega_i=n_i\omega$, where $n_i\in {N}$, $i=0,\ldots,n$, with $\omega$ defined as the greatest common factor of all the frequencies. Furthermore, assume that $n_i\not= n_j$ for $i\not =j$. Now, if we scale (\ref{eq-b2}) using $\tau = \omega t$, we get
\bqy
\nonumber
\omega\frac{d{\bf {\tilde{a}}}}{d\tau} = \frac{1}{2}{\bf Gd_1}(\tau)\left({\bf q}^T{\bf x}\right)^2+{\bf Gd_2}(\tau){\bf q}^T{\bf x}{\bf q}^T\left({\bf Ax}+\frac{1}{a_n}{\bf b}{\bf v}^T[\mathcal{S}(\tau)-\tilde{{\bf a}}]\right)
\eqy
with
\bqy
\nonumber	
{\bf d_1}(\tau)=\left[\begin{array}{c}
\sin(n_0\tau-\phi_0)+d_0\omega_0\cos(n_0\tau-\phi_0)\\
\sin(n_1\tau-\phi_1)+d_1\omega_1\cos(n_1\tau-\phi_1)\\
\vdots\\
\sin(n_n\tau-\phi_n)+d_n\omega_n\cos(n_n\tau-\phi_n)
\end{array}\right],\ 
{\bf d_2}(\tau)=\left[\begin{array}{c}
d_0\sin(n_0\tau-\phi_0)\\
d_1\sin(n_1\tau-\phi_1)\\
\vdots\\
d_n\sin(n_n\tau-\phi_n)
\end{array}\right]
\eqy 
and
\bqy
\mathcal{S}(\tau)=\left[\begin{array}{c}
c_0\sin n_0\tau\\
c_1\sin n_1\tau\\
\vdots\\
c_n\sin n_n\tau
\end{array}\right]
\eqy 
Defining $\varepsilon=1/\omega$, with $\omega\gg 1$, enables us to perform the averaging as follows
\bqy
\nonumber
{\bf \tilde{a}_{av}'}\triangleq\left(\frac{d{\bf {\tilde{a}}}}{d\tau}\right)_{\mbox{\scriptsize{av}}} &=& \varepsilon{\bf AVG}\left(\frac{1}{2}{\bf Gd_1}(\tau)\left({\bf q}^T{\bf x}\right)^2\right)\\
\nonumber
&&+\varepsilon{\bf AVG}\left({\bf Gd_2}(\tau){\bf q}^T{\bf x}{\bf q}^T\left[{\bf Ax}-\frac{1}{a_n}{\bf b}{\bf v}^T\tilde{{\bf a}}\right]\right)\\
&&+\varepsilon{\bf AVG}\left(\frac{1}{a_n}{\bf Gd_2}(\tau){\bf q}^T{\bf x}{\bf q}^T{\bf b}{\bf v}^T\mathcal{S}(\tau)\right)
\eqy
The first term can be written as
\bqy
\nonumber
{\bf AVG}\left(\frac{1}{2}{\bf Gd_1}(\tau)\left({\bf q}^T{\bf x}\right)^2\right)&=&
\frac{1}{2}{\bf G}\ {\bf AVG}\left({\bf d_1}(\tau)\left({\bf q}^T{\bf x}\right)^2\right)\\
\nonumber
&=&\frac{1}{2}{\bf G}\ {\bf AVG}\left({\bf d_1}(\tau)\right)\left({\bf q}^T{\bf x_{av}}\right)^2\\
\nonumber
&=& 0,
\eqy 
since ${\bf d_1}(\tau)$ is a zero mean periodic function over $T$. In a similar fashion, one can show that the second term also averages out to zero. Thus we are left with 
\bqy
{\bf \tilde{a}_{av}'} &=& \varepsilon{\bf AVG}\left(\frac{1}{a_n}{\bf Gd_2}(\tau){\bf q}^T{\bf x}{\bf q}^T{\bf b}{\bf v}^T\mathcal{S}(\tau)\right)
\eqy
Since ${\bf v}^T\mathcal{S}(\tau)$ is a scalar, we have ${\bf v}^T\mathcal{S}(\tau)=\mathcal{S}^T(\tau){\bf v}$. Thus we get
\bqy
\nonumber
{\bf \tilde{a}_{av}'} &=& \varepsilon{\bf AVG}\left(\frac{1}{a_n}{\bf G}{\bf d_2}(\tau)\mathcal{S}^T(\tau){\bf v}{\bf q}^T{\bf x}{\bf q}^T{\bf b}\right)\\
\label{eq-b8}
&=& \frac{\varepsilon}{a_n}{\bf G}\ {\bf AVG}\left({\bf d_2}(\tau)\mathcal{S}^T(\tau)\right){\bf v_{av}}{\bf q}^T{\bf x_{av}}{\bf q}^T{\bf b}
\eqy
However,
\bqy
\nonumber
&&{\bf d_2}(\tau)\mathcal{S}^T(\tau)=\\
\nonumber
&&\left[\begin{array}{ccc}
d_0c_0\sin(n_0\tau-\phi_0)\sin n_0\tau&\quad\cdots\quad&d_0c_n\sin(n_0\tau-\phi_0)\sin n_n\tau\\
\vdots&\quad\ddots\quad&\vdots\\
d_nc_0\sin(n_n\tau-\phi_n)\sin n_0\tau&\quad\cdots\quad&d_nc_n\sin(n_n\tau-\phi_n)\sin n_n\tau\\
\end{array}\right]
\eqy 
One can easily show that
\beq
\frac{1}{T}\int_0^{T=\frac{2\pi}{\omega}}\sin(n_i\tau-\phi_i)\sin n_j\tau\ d\tau=\left\{\begin{array}{ll}
0\quad &i\not=j\\
\frac{1}{2}\cos\phi_i\quad &i=j
\end{array}\right.
\eeq 
Thus
\bqy
&&{\bf AVG}\left({\bf d_2}(\tau)\mathcal{S}^T(\tau)\right)=
\frac{1}{2}\left[\begin{array}{cccc}
c_0d_0\cos\phi_0&0&\cdots&0\\
0&c_1d_1\cos\phi_1&\cdots&0\\
\vdots&\vdots&\ddots&\vdots\\
0&0&\cdots&c_nd_n\cos\phi_n
\end{array}\right]
\eqy 
Substituting this into (\ref{eq-b8}), and multiplying with the diagonal matrix $\bf G$ yields
\bqy
\nonumber
{\bf \tilde{a}_{av}'} = \frac{\varepsilon}{2a_n}\left[\begin{array}{cccc}
g_0c_0d_0\cos\phi_0&0&\cdots&0\\
0&g_1c_1d_1\cos\phi_1&\cdots&0\\
\vdots&\vdots&\ddots&\vdots\\
0&0&\cdots&g_nc_nd_n\cos\phi_n
\end{array}\right]{\bf v_{av}}{\bf q}^T{\bf x_{av}}{\bf q}^T{\bf b}
\eqy
Since, we have defined $\mathcal{C}=\mbox{diag}(g_id_ic_i\cos\phi_i)$, and since ${\bf q}^T{\bf b}=q_n$, the last equation simplifies to
\bqy
{\bf \tilde{a}_{av}'} = \varepsilon\frac{q_n}{2a_n}\mathcal{C}{\bf v_{av}}{\bf q}^T{\bf x_{av}}
\eqy
which is the same as (\ref{eq-systemII}).
% -------------------------------------------------------------------

\end{document}